\documentclass[proceedings]{dmtcs}

\usepackage[latin1]{inputenc}
\usepackage{subfigure}
\usepackage{pstricks}
\usepackage{fp}
\usepackage{multido}

\usepackage[round]{natbib}

\RCSdef$Revision: 1.3 $\endRCSdef
\rcsMajMin

\newtheorem{theo}{Theorem}

\newtheorem{lemm}[theo]{Lemma}

\newtheorem{defi}[theo]{Definition}

\newtheorem{rema}[theo]{Remark}

\newcommand{\pref}[1]{(\ref{#1})}

\def\cf{{\it cf. }}
\def\ie{{\it i.e. }}

\def\si{\sigma}
\def\ss{\sigma(a)^{-1}}

\def\bx{\bar x}

\def\ba{\bar a}

\newcommand{\flecheD}{\psline{->}(0,0)(.75,0)\psline(0,0)(1,0)}
\newcommand{\flecheG}{\psline{->}(0,0)(-.75,0)\psline(0,0)(-1,0)}
\newcommand{\flecheH}{\psline{->}(0,0)(0,.75)\psline(0,0)(0,1)}
\newcommand{\flecheB}{\psline{->}(0,0)(0,-.75)\psline(0,0)(0,-1)}

\newcounter{ISL}
\newcounter{ISH}
\newcommand{\IceGrid}[2]{%
\setcounter{ISL}{#1}\addtocounter{ISL}{-1}%
\setcounter{ISH}{#2}\addtocounter{ISH}{-1}%
\multido{\i=0+1}{#1}{\psline(\i,0)(\i,\theISH)}%
\multido{\i=0+1}{#2}{\psline(0,\i)(\theISL,\i)}%
}

\newcommand{\colD}[1]{%
\multido{\i=0+1}{#1}{\rput(0,\i){\flecheD}}%
}

\newcommand{\colG}[1]{%
\multido{\i=0+1}{#1}{\rput(1,\i){\flecheG}}%
}

\newcommand{\linH}[1]{%
\multido{\i=0+1}{#1}{\rput(\i,0){\flecheH}}%
}

\newcommand{\linB}[1]{%
\multido{\i=0+1}{#1}{\rput(\i,1){\flecheB}}%
}

\newcommand{\col}[1]{%
\multido{\i=0+1}{#1}{\rput(0,\i){\psline(0,0)(1,0)}}%
}

\newcommand{\lin}[1]{%
\multido{\i=0+1}{#1}{\rput(\i,0){\psline(0,0)(0,1)}}%
}

\newcommand{\halfCC}[5]{%
\degrees[360]%
\multido{\n=#3+#4}{#5}{\rput(#1,#2){\psarc(0,0){\n}{-90}{90}}}%
}

\newcommand{\tqCC}[5]{%
\degrees[360]%
\multido{\n=#3+#4}{#5}{\rput(#1,#2){\psarc(0,0){\n}{-90}{180}}}%
}

\newcommand{\IceSquare}[1]{%
\rput(1,1){\IceGrid{#1}{#1}}%
\rput(0,1){\colD{#1}}%
\rput(1,0){\linB{#1}}%
\rput(#1,1){\colG{#1}}%
\rput(1,#1){\linH{#1}}%
}

\newcounter{ISdoublesize}
\newcommand{\HTIceEven}[1]{%
\setcounter{ISdoublesize}{#1}\addtocounter{ISdoublesize}{#1}%
\rput(1,1){\IceGrid{#1}{\theISdoublesize}}%
\rput(0,1){\colD{\theISdoublesize}}%
\rput(1,0){\linB{#1}}%
\rput(1,\theISdoublesize){\linH{#1}}%
\rput(#1,#1){\rput(0,.5){\halfCC{0}{0}{.5}{1.0}{#1}}}%
\rput(#1,#1){\psline[linestyle=dotted](.25,.5)(.75,.5)}%
}

\newcommand{\HTIceOdd}[1]{%
\setcounter{ISdoublesize}{#1}\addtocounter{ISdoublesize}{#1}%
\addtocounter{ISdoublesize}{1}%
\rput(1,1){\IceGrid{#1}{\theISdoublesize}}%
\rput(0,1){\colD{\theISdoublesize}}%
\rput(1,0){\linB{#1}}%
\rput(1,\theISdoublesize){\linH{#1}}%
\rput(#1,0){\psline(1,1)(1,#1)\rput(1,1){\flecheB}}%
\psarc(#1,#1){1}{0}{90}%
\rput(#1,#1){\halfCC{1}{1}{1}{1}{#1}}%
\rput(#1,1){\col{#1}}%
\rput(#1,#1){\rput(0,2){\col{#1}}}%
\rput(#1,#1){\psline[linestyle=dotted](.25,.25)(1.25,1.25)}
}

\newcounter{QTsize}
\newcommand{\QTIce}[1]{%
\setcounter{QTsize}{#1}\addtocounter{QTsize}{-1}%
\rput(1,1){\IceGrid{\theQTsize}{\theQTsize}}%
\rput(0,1){\colD{#1}}%
\rput(1,0){\linB{#1}}%
\rput(\theQTsize,1){\col{\theQTsize}}%
\rput(1,\theQTsize){\lin{\theQTsize}}%
\rput(\theQTsize,\theQTsize){\psarc(0,0){1}{0}{90}}%
\psline(1,#1)(\theQTsize,#1)\psline(#1,1)(#1,\theQTsize)%
\rput(#1,#1){\tqCC{0}{0}{1}{1}{\theQTsize}}%
\rput(#1,#1){\SpecialCoor\multido{\i=1+1}{\theQTsize}{\rput(\i;45){\psdots[dotstyle=*](0,0)}}}%
\rput(\theQTsize,\theQTsize){\SpecialCoor\psdots[dotstyle=*](1;45)\psline[linestyle=dotted](.5;45)(1.5;45)}%
}

\newcommand{\qQTIce}[1]{%
\setcounter{QTsize}{#1}\addtocounter{QTsize}{-1}%
\rput(1,1){\IceGrid{\theQTsize}{\theQTsize}}%
\rput(0,1){\colD{#1}}%
\rput(1,0){\linB{#1}}%
\rput(\theQTsize,1){\col{\theQTsize}}%
\rput(1,\theQTsize){\lin{\theQTsize}}%
\rput(\theQTsize,\theQTsize){\psarc(0,0){1}{0}{90}}%
\psline(1,#1)(\theQTsize,#1)\psline(#1,1)(#1,\theQTsize)%
\rput(#1,#1){\tqCC{0}{0}{1}{1}{\theQTsize}}%
\rput(#1,#1){\SpecialCoor\multido{\i=1+1}{\theQTsize}{\rput(\i;45){\psdots[dotstyle=*](0,0)}}}%
\rput(\theQTsize,\theQTsize){\SpecialCoor\psline[linestyle=dotted](.5;45)(1.5;45)}%
}

\newcommand{\Hcrossing}{%
\psbezier(0,0)(.4,0)(.6,1)(1,1)%
\psbezier(0,1)(.4,1)(.6,0)(1,0)%
}

\newcommand{\convcorner}{%
\begin{pspicture}(0.2,0.2)\psset{linewidth=.4pt,arrowsize=2.5pt}
\psline{->}(0,.2)(.2,.2)
\psline{->}(.2,0)(.2,.2)
\end{pspicture}%
}

\newcommand{\divcorner}{%
\begin{pspicture}(.2,.2)\psset{linewidth=.4pt,arrowsize=2.5pt}
\psline{->}(.2,.2)(.2,0)
\psline{->}(.2,.2)(0,.2)
\end{pspicture}%
}

\newcommand{\upleft}{%
\begin{pspicture}(.2,.2)\psset{linewidth=.4pt,arrowsize=2.5pt}
\psline(.2,0)(.2,.15)\psarc(.15,.15){.05}{0}{90}\psline{->}(.15,.2)(0,.2)%
\end{pspicture}%
}

\newcommand{\downright}{%
\begin{pspicture}(.2,.2)\psset{linewidth=.4pt,arrowsize=2.5pt}
\psline(0,.2)(.15,.2)\psarc(.15,.15){.05}{0}{90}\psline{->}(.2,.15)(.2,0)%
\end{pspicture}
}

\newcommand{\upc}{%
\begin{pspicture}(.2,.2)\psset{linewidth=.4pt,arrowsize=2.5pt}%
\psline(0,0)(.1,0)\psarc(.1,.1){.1}{270}{90}\psline{->}(.1,.2)(0,.2)%
\end{pspicture}%
}

\newcommand{\downc}{%
\begin{pspicture}(.2,.2)\psset{linewidth=.4pt,arrowsize=2.5pt}%
\psline(0,.2)(.1,.2)\psarc(.1,.1){.1}{270}{90}\psline{->}(.1,0)(0,0)%
\end{pspicture}%
}

\newcommand{\ZqConv}{%
Z_{\textsc{QT}}^{\convcorner}%
}

\newcommand{\ZqDiv}{%
Z_{\textsc{QT}}^{\divcorner}%
}

\newcommand{\ZqUpleft}{%
Z_{\textsc{QT}}^{\upleft}%
}

\newcommand{\ZqDownright}{%
Z_{\textsc{QT}}^{\downright}%
}

\newcommand{\ZhUp}{%
Z_{\textsc{HT}}^{\upc}%
}

\newcommand{\ZhDown}{%
Z_{\textsc{HT}}^{\downc}%
}

\newcommand{\ZhUpleft}{%
Z_{\textsc{HT}}^{\upleft}%
}

\newcommand{\ZhDownright}{%
Z_{\textsc{HT}}^{\downright}%
}

\author{Jean-Christophe Aval, Philippe Duchon\thanks{Both authors are supported by the ANR project MARS (BLAN06-2$\_$0193)} 
}

\address{LaBRI, Universit\'e Bordeaux 1, CNRS\\ 351 cours
 de la Lib\'eration, 33405 Talence cedex, FRANCE}

\title[Quarter-turn Symmetric Alternating Sign Matrices]{Enumeration of alternating sign matrices of even size (quasi)-invariant under a quarter-turn rotation}

\date{\today}

\begin{document}
\maketitle

\begin{abstract}
\noindent {\bf \normalsize Abstract.}
The aim of this work is to enumerate alternating sign matrices (ASM) that are quasi-invariant under a quarter-turn.
The enumeration formula (conjectured by Duchon) involves, as a product of three terms, the number of unrestricted ASM's and the number of half-turn symmetric ASM's.

\noindent {\bf \normalsize R\'esum\'e.}
L'objet de ce travail est d'\'enum\'erer les matrices \`a signes alternants (ASM) quasi-invariantes par rotation d'un quart-de-tour. La formule d'\'enum\'eration, conjectur\'ee par Duchon, fait appara\^itre trois facteurs, comprenant le nombre d'ASM quelconques et le nombre d'ASM invariantes par demi-tour.

\end{abstract}

\section{Introduction}\label{sec:intro}

An {\em alternating sign matrix} is a square matrix with entries in $\{-1,0,1\}$ and such that in any row and column: the non-zero entries alternate in sign, and their sum is equal to $1$. Their enumeration formula was conjectured by Mills, Robbins and Rumsey \cite{MRR}, and proved by Zeilberger \cite{zeil}, and almost simultaneously by Kuperberg \cite{kup1}. Kuperberg used a bijection between the ASM's and the states of a statistical square-ice model, for which he studied and computed the partition function. He also used this method in \cite{kup} to obtain many enumeration or equinumeration results for various symmetry classes of ASM's, most of them having been conjectured by Robbins \cite{robbins}. Among these results can be found the following remarkable one.

\begin{theo}\label{theo:kup}
{\em (Kuperberg).}
The number $A_{\textsc{QT}}(4N)$ of ASM's of size $4N$ invariant under a quarter-turn (QTASM's) is related to the number $A(N)$ of (unrestricted) ASM's of size $N$ and to the number $A_{\textsc{HT}}(2N)$ of ASM's of size $2N$ invariant under a half-turn (HTASM's) by the formula:
\begin{equation}\label{eq:kup}
A_{\textsc{QT}}(4N)=A_{\textsc{HT}}(2N) A(N)^2.
\end{equation}
\end{theo}

More recently, Razumov and Stroganov \cite{RS} applied Kuperberg's strategy to settle the following result relative to QTASM's of odd size, also conjectured by Robbins \cite{robbins} .

\begin{theo}\label{theo:RS}
{\em (Razumov, Stroganov).}
The numbers of QTASM's of odd size are given by the following formulas, where $A_{\textsc{HT}}(2N+1)$ is the number of HTASM's of size $2N+1$:
  \begin{eqnarray}
    A_{\textsc{QT}}(4N-1) & = & A_{\textsc{HT}}(2N-1) A(N)^2\label{eq:QT_m1}\\
    A_{\textsc{QT}}(4N+1) & = & A_{\textsc{HT}}(2N+1) A(N)^2\label{eq:QT_p1}.
  \end{eqnarray}
\end{theo}

It is easy to observe (and will be proved in Section \ref{sec:qQTASM}) that the set of QTASM's of size $4N+2$ is empty. But, by slightly relaxing the symmetry condition at the center of the matrix, Duchon introduced in \cite{duchon} the notion of ASM's quasi-invariant under a quarter turn (the definition will be given in Section \ref{sec:qQTASM}) whose class is non-empty in size $4N+2$. Moreover, he conjectured for these qQTASM's an enumeration formula that perfectly completes the three previous enumeration results on QTASM. It is the aim of this paper to establish this formula.

\begin{theo}\label{theo:form}
The number $A_{\textsc{QT}}(4N+2)$ of qQTASM of size $4N+2$ is given by:
\begin{equation}\label{eq:form}
A_{\textsc{QT}}(4N+2)=A_{\textsc{HT}}(2N+1) A(N)A(N+1).
\end{equation}
\end{theo}

This paper is organized as follows: in Section \ref{sec:qQTASM}, we define qQTASM's; in Section \ref{sec:icemodel}, we recall the definitions of square ice models, precise the parameters and the partition functions that we shall study, and give the formula corresponding to equation \pref{eq:form} at the level of partition functions; Section \ref{sec:proofs} is devoted to the proofs.

\section{ASM's quasi-invariant under a quarter-turn}\label{sec:qQTASM}

The class of ASM's invariant under a rotation by a quarter-turn (QTASM) is non-empty in size $4N-1$, $4N$, and $4N+1$. But this is not the case in size $4N+2$.

\begin{lemm} 
There is no QTASM of size $4N+2$.
\end{lemm}
\proof
Let us suppose that $M$ is a QTASM of even size $2L$. Now we use the fact that the size of an ASM is given by the sum of its entries, and the symmetry of $M$ to write:
\begin{equation}
2L=\sum_{1\le i,j\le2L}M_{i,j}=4\times\sum_{1\le i,j\le L}M_{i,j}
\end{equation}
which implies that the size of $M$ has to be a multiple of $4$.
\endproof

Duchon introduced in \cite{duchon} a notion of ASM's quasi-invariant under a quarter-turn, by slightly relaxing the symmetry condition at the center of the matrix. The definition is more simple when considering the height matrix associated to the ASM, but can also be given directly.

\begin{defi}
An ASM $M$ of size $4N+2$ is said to be {\em quasi-invariant under a quarter-turn} (qQTASM) if its entries satisfy the quarter-turn symmetry
\begin{equation}
M_{4N+2-j+1,4N+2-i+1}=M_{i,j}
\end{equation}
except for the four central entries $(M_{2N,2N},M_{2N,2N+1},M_{2N+1,2N},M_{2N+1,2N+1})$ that have to be either $(0,-1,-1,0)$ or $(1,0,0,1)$.
\end{defi}

We give below two examples of qQTASM's of size $6$, with the two possible patterns at the center.
$$\left(
\begin{array}{cccccc}
0&0&0&1&0&0\cr
0&0&1&0&0&0\cr
1&0&0&-1&1&0\cr
0&1&-1&0&0&1\cr
0&0&0&1&0&0\cr
0&0&1&0&0&0\cr
\end{array}
\right)
\ \ \ \ \ \ \ \ \ \ \ \ \ \ \ \ 
\left(
\begin{array}{cccccc}
0&0&1&0&0&0\cr
0&1&-1&0&1&0\cr
0&0&1&0&-1&1\cr
1&-1&0&1&0&0\cr
0&1&0&-1&1&0\cr
0&0&0&1&0&0\cr
\end{array}
\right)$$

In the next section, we associate square ice models to ASM's with various types of symmetry.

\section{Square ice models and partition functions}\label{sec:icemodel}

\subsection{Notations}

Using Kuperberg's method we introduce square ice models associated to ASM's, HTASM's and (q)QTASM's. We recall here the main definitions and refer to \cite{kup} for details and many examples. 

Let $a\in\mathbb{C}$ be a global parameter.
For any complex number $x$ different from zero, we denote $\overline{x}=1/x$, and we define:
\begin{equation}
  \sigma(x) =  x - \overline{x}.
\end{equation}

If $G$ is a tetravalent graph, an {\em ice state} of $G$ is an orientation of the edges such that every tetravalent vertex has exactly two incoming and two outgoing edges.

\begin{figure}[htbp]
  \begin{center}
\psset{unit=0.5cm}
\begin{pspicture}[.5](3,3)
\rput(1.5,1.5){
   \psline(-1,0)(1,0)
   \psline(0,-1)(0,1)
   \psdot[dotsize=.2]
   \rput(-.5,-.5){$x$}
  }
\end{pspicture}
$=$
    \begin{pspicture}[.5](19,5)
\rput(2,2.5){%
    \rput(-1,0){\flecheD}%
    \rput(1,0){\flecheG}%
    \rput(0,0){\flecheH}%
    \rput(0,0){\flecheB}%
    \rput[t](0,-1.2){$\sigma(a^2)$}%
    \rput[t](0,-2.2){$1$}%
}
\rput(5,2.5){%
    \rput(0,0){\flecheG}%
    \rput(0,0){\flecheD}%
    \rput(0,-1){\flecheH}%
    \rput(0,1){\flecheB}%
    \rput[t](0,-1.2){$\sigma(a^2)$}%
    \rput[t](0,-2.2){$-1$}%
}
\rput(8,2.5){%
    \rput(0,0)\flecheG
    \rput(0,0)\flecheH
    \rput(1,0)\flecheG
    \rput(0,-1)\flecheH
    \rput[t](0,-1.2){$\sigma(ax)$}
    \rput[t](0,-2.2){$0$}%
}
\rput(11,2.5){
    \rput(0,0)\flecheD
    \rput(0,0)\flecheB
    \rput(-1,0)\flecheD
    \rput(0,1)\flecheB
    \rput[t](0,-1.2){$\sigma(ax)$}
    \rput[t](0,-2.2){$0$}%
}
\rput(14,2.5){
    \rput(-1,0)\flecheD
    \rput(0,0)\flecheD
    \rput(0,-1)\flecheH
    \rput(0,0)\flecheH
    \rput[t](0,-1.2){$\sigma(a\overline{x})$}
    \rput[t](0,-2.2){$0$}%
}
\rput(17,2.5){
    \rput(0,0)\flecheG
    \rput(1,0)\flecheG
    \rput(0,0)\flecheB
    \rput(0,1)\flecheB
    \rput[t](0,-1.2){$\sigma(a\overline{x})$}
    \rput[t](0,-2.2){$0$}%
}
    \end{pspicture}
  \end{center}
\caption{The 6 possible orientations, their associated weights and the corresponding entries in ASM's}
\label{fig:poids_6V}
\end{figure}
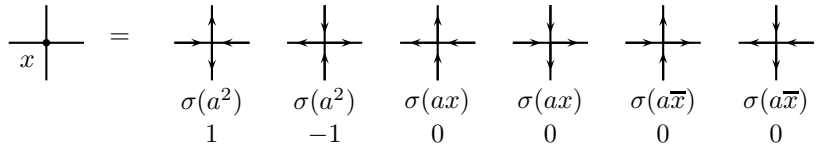

A parameter $x\neq 0$ is assigned to any tetravalent vertex of the graph $G$. 
Then this vertex gets a weight, which depends on its orientations, as shown on 
Figure~\ref{fig:poids_6V}.

It is sometimes easier to assign parameters, not to each vertex of the graph, but to the lines that compose the graph. In this case, the weight of a vertex is defined as:

\begin{displaymath}
\psset{unit=.5cm}
  \begin{pspicture}[.45](2,2)
\psline(1,0)(1,2)\psline(0,1)(2,1)
\rput[r](-.2,1){$x$}\rput[t](1,-.2){$y$}
  \end{pspicture}\ =\ 
  \begin{pspicture}[.45](2,2)
\psline(0,1)(2,1)\psline(1,0)(1,2)
\rput(.25,.25){$x\overline{y}$}
  \end{pspicture}
\end{displaymath}

When this convention is used, a parameter explicitly written at a vertex replaces the quotient of the parameters of the lines.

We will put a dotted line to mean that the parameter of a line is different on the two sides of the dotted line.
We will also use divalent vertices, and in this case the two edges have to be both in, or both out, and the corresponding weight is $1$:
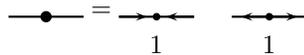
\begin{displaymath}
  \psset{unit=.5cm}
  \begin{pspicture}[.5](2,1)
    \psline(0,.5)(2,.5)\rput(1,.5){\psdot[dotsize=.3]}
  \end{pspicture} = 
  \begin{pspicture}[.5](5,1)
    \rput(0,.5)\flecheD
    \rput(2,.5)\flecheG
    \rput[t](1,0){$1$}
    \rput(4,.5)\flecheG
    \rput(4,.5)\flecheD
    \rput[t](4,0){$1$}
    \rput(1,.5){\psdot[dotsize=.2]}
    \rput(4,.5){\psdot[dotsize=.2]}
  \end{pspicture}
\end{displaymath}

The partition function of a given ice model is then defined as the sum over all its states of the product of the weights of the vertices.

To simplify notations, we will denote by $X_{N}$ the vector of variables
$(x_1,\dots, x_N)$. We use the notation
$X\backslash x$ to denote the vector $X$ without the variable $x$.

\subsection{Partition functions for classes of ASM's}

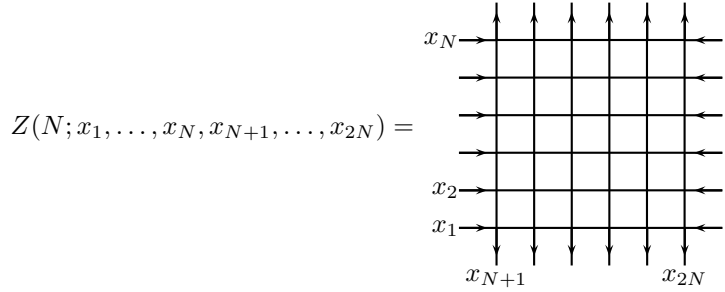
\begin{figure}[htbp]
  \begin{displaymath}
    Z(N;x_1,\dots, x_{N}, x_{N+1},\dots, x_{2N}) =
    \psset{unit=.5cm}
    \begin{pspicture}[.5](7,7)
      \rput(1,0){\IceSquare{6}}
      \rput[r](1,1){$x_{1}$}\rput[r](1,2){$x_{2}$}\rput[r](1,6){$x_{N}$}
      \rput[t](2,-.1){$x_{N+1}$}\rput[t](3,-.1){}\rput[t](7,-.1){$x_{2N}$}
    \end{pspicture}
  \end{displaymath}
\caption{Partition function for ASM's of size $N$}
\label{fig:Z}
\end{figure}

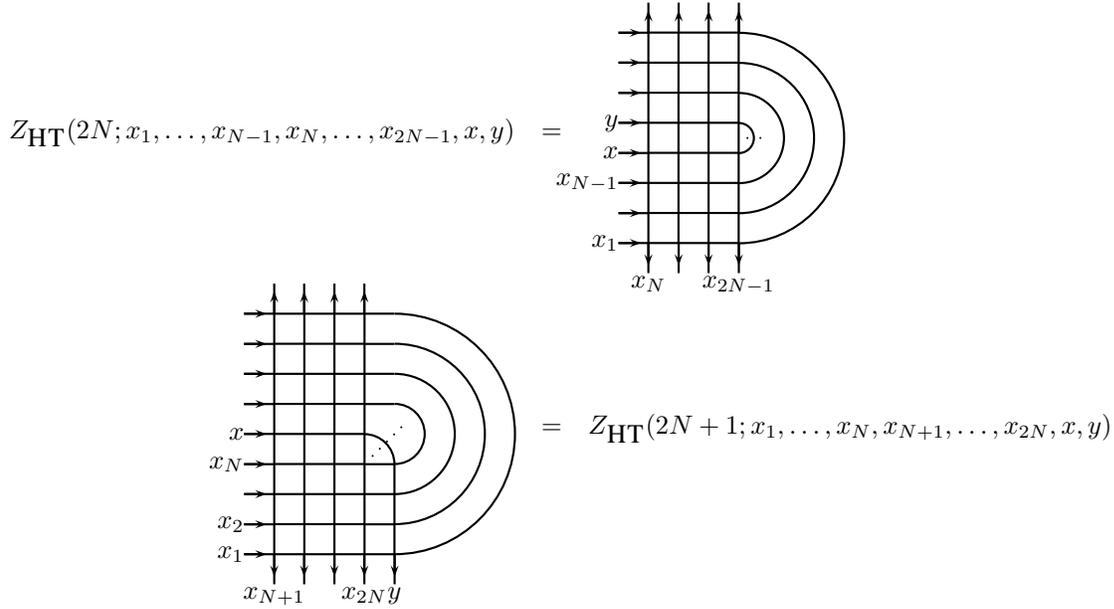
\begin{figure}[htbp]
  \begin{eqnarray*}
    Z_{\textsc{HT}}(2N; x_1,\dots,x_{N-1},x_{N},\dots,x_{2N-1},x,y)&  =& 
\psset{unit=.4cm}
    \begin{pspicture}[.5](9,9)
      \rput(1,0){\HTIceEven{4}}
      \rput[r](1,1){$x_{1}$}\rput[r](1,3){$x_{N-1}$}\rput[r](1,5){$y$}
      \rput[r](1,4){$x$}
      \rput[t](2,-.1){$x_{N}$}\rput[t](3,-.1){}\rput[t](5,-.1){$x_{2N-1}$}
    \end{pspicture}\\ 
\psset{unit=.4cm}
\begin{pspicture}[.5](10,10)
  \rput(1,0){\HTIceOdd{4}}
  \rput[r](1,1){$x_{1}$}\rput[r](1,2){$x_{2}$}\rput[r](1,4){$x_{N}$}
  \rput[r](1,5){$x$}
  \rput[t](2,-.1){$x_{N+1}$}\rput[t](3,-.1){}\rput[t](5,-.1){$x_{2N}$}
  \rput[t](6,-.1){$y$}
\end{pspicture} & = & 
Z_{\textsc{HT}}(2N+1;x_1,\dots,x_{N},x_{N+1},\dots,x_{2N},x,y) 
  \end{eqnarray*}
\caption{Partition functions for HTASM's}
\label{fig:ZHT}
\end{figure}

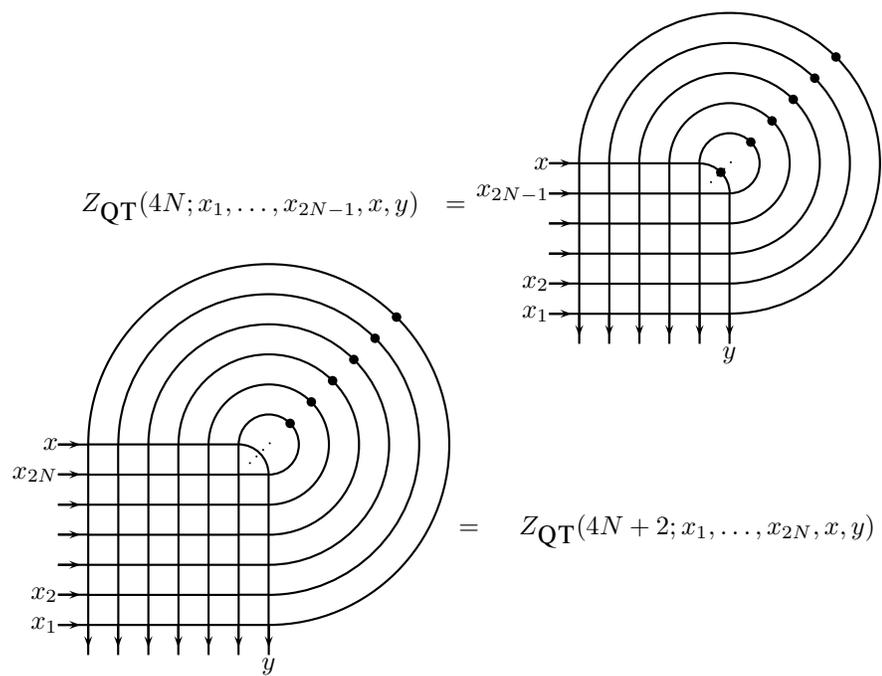
\begin{figure}[htbp]
\begin{eqnarray*}
  Z_{\textsc{QT}}(4N;x_1,\dots,x_{2N-1},x,y) & =\ \  \ \  & 
\psset{unit=.4cm}
  \begin{pspicture}[.4](11,11)
    \rput(1,0){\QTIce{6}}
    \rput[r](1,1){$x_1$}\rput[r](1,2){$x_2$}\rput[r](1,5){$x_{2N-1}$}
    \rput[r](1,6){$x$}\rput[t](7,-.1){$y$}
  \end{pspicture}\\
\psset{unit=.4cm}
  \begin{pspicture}[.4](12,10)
    \rput(0,0){\qQTIce{7}}
    \rput[r](0,1){$x_1$}\rput[r](0,2){$x_2$}\rput[r](0,6){$x_{2N}$}
    \rput[r](0,7){$x$}\rput[t](7,-.1){$y$}
  \end{pspicture} 
  & = & Z_{\textsc{QT}}(4N+2; x_1,\dots, x_{2N},x,y)
\end{eqnarray*}
\caption{Partition functions for (q)QTASM of even size}
\label{fig:ZQT}
\end{figure}

We give in Figures \ref{fig:Z}, \ref{fig:ZHT}, and \ref{fig:ZQT} 
the ice models corresponding to the classes of ASM's that we shall study, and their partition functions. The bijection between (unrestricted) ASM's and states of the square ice model with ``domain wall boundary'' is now well-known (\cf \cite{kup}), and the bijections for the other symmetry classes may be easily checked in the same way. The correspondence between orientations of the ice model and entries of ASM's is given in Figure \ref{fig:poids_6V}.

With these notations, Theorem \ref{theo:form} will be a consequence of the following one which addresses the concerned partition functions.

\begin{theo}\label{theo:main}
When $a=\omega_6=\exp(i\pi/3)$, one has for $N\ge 1$:
\begin{equation}\label{eq:main1}
Z_{\textsc{QT}}(4N;X_{2N-1},x,y)=\ss Z_{\textsc{HT}}(2N;X_{2N-1},x,y)Z(N;X_{2N-1},x)Z(N;X_{2N-1},y)
\end{equation}
and
\begin{equation}\label{eq:main2}
Z_{\textsc{QT}}(4N+2;X_{2N},x,y)=\ss Z_{\textsc{HT}}(2N+1;X_{2N},x,y)Z(N;X_{2N})Z(N+1;X_{2N},x,y).
\end{equation}
\end{theo}

Equation \pref{eq:main2} is new; equation \pref{eq:main1} is due to Kuperberg \cite{kup} for the case $x=y$.
To see that Theorem \ref{theo:main} implies Theorem \ref{theo:form} (and Theorem \ref{theo:kup}), we just have to observe that when $a=\omega_6$ and all the variables are set to $1$, then the weight at each vertex is $\si(a)=\si(a^2)$ thus the partition function reduces (up to multiplication by $\si(a)^{\rm number\ of\ vertices}$) to the number of states.

\section{Proofs}\label{sec:proofs}

In this extended abstract, we shall only give the main ideas of the proofs. Most of them are greatly inspired from \cite{kup}.
To prove Theorem \ref{theo:main}, the method is to identify both sides of equations \pref{eq:main1} and \pref{eq:main2} as Laurent polynomials, and to produce as many specializations of the variables that verify the equalities, as needed to imply these equations in full generality.

\subsection{Laurent polynomials}

Since the weight of any vertex is a Laurent polynomial in the variables $x_i$, $x$ and $y$, the partition functions are Laurent polynomials in these variables. Moreover they are centered Laurent polynomials, \ie their lowest degree is the negative of their highest degree (called the half-width of the polynomial). In order to divide by two the number of non-zero coefficients (hence the number of required specializations) in $x$, we shall deal with Laurent polynomials of given parity in this variable. To do so, we group together the states with a given orientation (indicated as subscripts in the following notations) at the edge where the parameters $x$ and $y$ meet.

So let us consider the partition functions $\ZqConv(4N;X_{2N-1},x,y)$ and
$\ZqDiv(4N;X_{2N-1},x,y)$, respectively odd and even parts of
$Z_{\textsc{QT}}(4N;X_{2N-1},x,y)$ in $x$;
$\ZqUpleft(4N+2;X_{2N},x,y)$ and
$\ZqDownright(4N+2;X_{2N},x,y)$, respectively odd and even parts of
$Z_{\textsc{QT}}(4N+2;X_{2N},x,y)$in $x$;
$\ZhUp(2N;X_{2N-1},x,y)$ and
$\ZhDown(2N;X_{2N-1},x,y)$, respectively parts with the parity of $N$ and of $N-1$ of
$Z_{\textsc{HT}}(2N;X_{2N-1},x,y)$ in $x$; and
$\ZhUpleft(2N+1;X_{2N},x,y)$ and
$\ZhDownright(2N+1;X_{2N},x,y)$, respectively parts with the parity of $N-1$ and of $N$ of
$Z_{\textsc{HT}}(2N+1;X_{2N},x,y)$ in $x$.

With these notations, the equations \pref{eq:main1} and \pref{eq:main2}
are equivalent to the following:
\begin{eqnarray}
 \!\!\!\!\!\! \sigma(a)\ZqConv(4N;X_{2N-1},x,y) &\!\! =\!\! &
  \ZhUp(2N;X_{2N-1},x,y) Z(N;X_{2N-1},x)
  Z(N;X_{2N-1},y), \label{eq:main1.1}\\
 \!\!\!\!\!\! \sigma(a)\ZqDiv(4N;X_{2N-1},x,y) &\!\! =\!\! &
  \ZhDown(2N;X_{2N-1},x,y) Z(N;X_{2N-1},x)
  Z(N;X_{2N-1},y), \label{eq:main1.2}\\
 \!\!\!\!\!\!  \sigma(a)\ZqDownright(4N+2;X_{2N},x,y) &\!\! =\!\! &
  \ZhDownright(2N+1;X_{2N},x,y)
  Z(N+1;X_{2N},x,y) Z(N;X_{2N}), \label{eq:main2.1}\\
 \!\!\!\!\!\! \sigma(a)\ZqUpleft(4N+2;X_{2N},x,y) &\!\! =\!\! &
  \ZhUpleft(2N+1;X_{2N},x,y)
  Z(N+1;X_{2N},x,y) Z(N;X_{2N}). \label{eq:main2.2}
\end{eqnarray}

\begin{lemm}\label{lemm:Laurent}
Both left-hand side and right-hand side of equations (\ref{eq:main1.1}-\ref{eq:main2.2}) are centered Laurent polynomials 
in the variable $x$,
odd or even, of respective half-widths $2N-1$, $2N-2$, $2N$, and
$2N-1$. Thus to prove each of these identities we have to exhibit specializations of $x$ 
for which the equality is true, and in number strictly exceeding the half-width.
\end{lemm}
\proof
To compute the half-width of these partition functions, just count the number of vertices in the ice models, and take note that non-zero entries of the ASM (\ie the first two orientations of Figure  \ref{fig:poids_6V}) give constant weight $\si(a^2)$. Also, a line whose orientation changes between endpoints must have an odd (hence non-zero) number of these $\pm 1$ entries.
\endproof

\subsection{Symmetries}

To produce many specializations from one, we shall use symmetry properties of the partition functions. The crucial tool to prove this is the Yang-Baxter equation that we recall below.

\begin{lemm}\label{lemm:YB}{\em [Yang-Baxter equation]}
  If $xyz=\overline{a}$, then
  \begin{equation}
    \begin{pspicture}[.5](2,2)
      \SpecialCoor
      \rput(1,1){
	\psarc(-1.732,0){2}{330}{30}
        \psarc(.866,1.5){2}{210}{270}
        \psarc(.866,-1.5){2}{90}{150}
	\rput(1;60){$x$}
        \rput(1;180){$y$}
        \rput(1;300){$z$}
      }
    \end{pspicture}
     = 
     \begin{pspicture}[.5](2,2)
      \rput(1,1){
	\SpecialCoor
	\psarc(1.732,0){2}{150}{210}
        \psarc(-.866,1.5){2}{270}{330}
        \psarc(-.866,-1.5){2}{30}{90}
        \rput(1;240){$x$}
        \rput(1;0){$y$}
        \rput(1;120){$z$}
      }
     \end{pspicture}.
\label{eq:YB}
  \end{equation}
\end{lemm}

The following lemma gives a (now classical) example of use of the Yang-Baxter equation.

\begin{lemm}\label{lemm:echange_lignes}
  \begin{equation}
    \psset{unit=.5cm}
    \begin{pspicture}[.5](5.5,2)
      \rput[r](1,.5){$x$}\rput[r](1,1.5){$y$}
      \rput(1,.5){\flecheD}\rput(1,1.5){\flecheD}
      \psline(2,0)(2,2)\psline(3,0)(3,2)\psline(4.5,0)(4.5,2)
      \psline(2,.5)(3.5,.5)\psline(2,1.5)(3.5,1.5)
      \psline(4,.5)(4.5,.5)\psline(4,1.5)(4.5,1.5)
      \rput(3.75,1){$\dots$}
      \rput(5.5,.5){\flecheG}
      \rput(5.5,1.5){\flecheG}
    \end{pspicture} = 
    \begin{pspicture}[.5](5.5,2)
      \rput[r](1,.5){$y$}\rput[r](1,1.5){$x$}
      \rput(1,.5){\flecheD}\rput(1,1.5){\flecheD}
      \psline(2,0)(2,2)\psline(3,0)(3,2)\psline(4.5,0)(4.5,2)
      \psline(2,.5)(3.5,.5)\psline(2,1.5)(3.5,1.5)
      \psline(4,.5)(4.5,.5)\psline(4,1.5)(4.5,1.5)
      \rput(3.75,1){$\dots$}
      \rput(5.5,.5){\flecheG}
      \rput(5.5,1.5){\flecheG}      
    \end{pspicture}.
    \label{eq:symetrie_Z}
  \end{equation}
\end{lemm}

\proof
We multiply the left-hand side by $\sigma(a\overline{z})$, with
  $z=\overline{a}x\overline{y}$. We get
  \begin{eqnarray*} \psset{unit=.5cm}
    \sigma(a\overline{z})
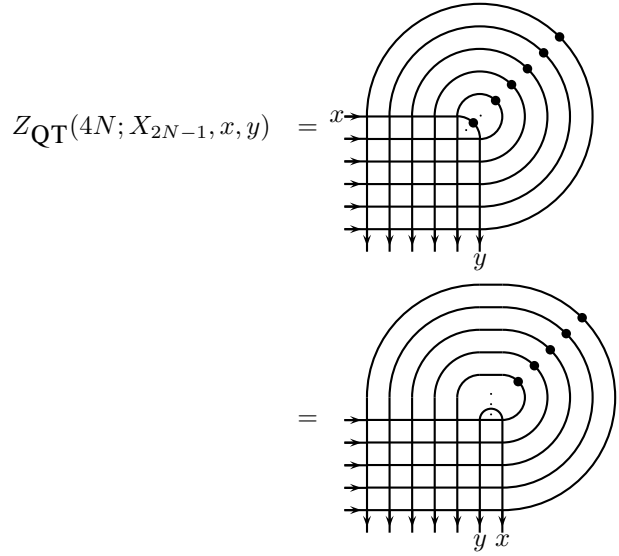
    \begin{pspicture}[.5](5.5,2)
      \rput[r](1,.5){$x$}\rput[r](1,1.5){$y$}
      \rput(1,.5){\flecheD}\rput(1,1.5){\flecheD}
      \psline(2,0)(2,2)\psline(3,0)(3,2)\psline(4.5,0)(4.5,2)
      \psline(2,.5)(3.5,.5)\psline(2,1.5)(3.5,1.5)
      \psline(4,.5)(4.5,.5)\psline(4,1.5)(4.5,1.5)
      \rput(3.75,1){$\dots$}
      \rput(5.5,.5){\flecheG}
      \rput(5.5,1.5){\flecheG}
    \end{pspicture} & = & \psset{unit=.5cm}
    \begin{pspicture}[.5](6.5,2)
      \rput[r](1,.5){$y$}\rput[r](1,1.5){$x$}
      \rput(1,.5){\flecheD}\rput(1,1.5){\flecheD}
      \rput(2,.5){\Hcrossing}\rput[r](2.4,1){$z$}
      \rput(1,0){
	\psline(2,0)(2,2)\psline(3,0)(3,2)\psline(4.5,0)(4.5,2)
	\psline(2,.5)(3.5,.5)\psline(2,1.5)(3.5,1.5)
	\psline(4,.5)(4.5,.5)\psline(4,1.5)(4.5,1.5)
	\rput(3.75,1){$\dots$}
	\rput(5.5,.5){\flecheG}
	\rput(5.5,1.5){\flecheG}}
    \end{pspicture} \\ 
      & = & 
    \psset{unit=.5cm}
    \begin{pspicture}[.5](5.5,2)
      \rput[r](1,.5){$y$}\rput[r](1,1.5){$x$}
      \rput(1,.5){\flecheD}\rput(1,1.5){\flecheD}
      \psline(2,0)(2,2)
      \rput(2,.5){\Hcrossing}\rput[r](2.4,1){$z$}
      \psline(3,0)(3,2)\psline(4.5,0)(4.5,2)
      \psline(3,.5)(3.5,.5)\psline(3,1.5)(3.5,1.5)
      \rput(3.75,1){$\dots$}
      \psline(4,.5)(4.5,.5)\psline(4,1.5)(4.5,1.5)
      \rput(5.5,.5){\flecheG}
      \rput(5.5,1.5){\flecheG}
    \end{pspicture} \\
      & = & 
    \psset{unit=.5cm}
    \begin{pspicture}[.5](6.5,2)
      \rput[r](1,.5){$y$}\rput[r](1,1.5){$x$}
      \rput(1,.5){\flecheD}\rput(1,1.5){\flecheD}
      \psline(2,0)(2,2)
      \psline(3,0)(3,2)
      \psline(2,.5)(3.5,.5)\psline(2,1.5)(3.5,1.5)
      \rput(3.75,1){$\dots$}
      \psline(4,.5)(4.5,.5)\psline(4,1.5)(4.5,1.5)
      \rput(4.5,.5){\Hcrossing}\rput[r](4.9,1){$z$}
      \psline(5.5,0)(5.5,2)
      \rput(6.5,.5){\flecheG}
      \rput(6.5,1.5){\flecheG}
    \end{pspicture} \\
      & = & 
    \psset{unit=.5cm}
    \begin{pspicture}[.5](6.5,2)
      \rput[r](1,.5){$y$}\rput[r](1,1.5){$x$}
      \rput(1,.5){\flecheD}\rput(1,1.5){\flecheD}
      \psline(2,0)(2,2)\psline(3,0)(3,2)\psline(4.5,0)(4.5,2)
      \psline(2,.5)(3.5,.5)\psline(2,1.5)(3.5,1.5)
      \psline(4,.5)(4.5,.5)\psline(4,1.5)(4.5,1.5)
      \rput(3.75,1){$\dots$}
      \rput(4.5,.5){\Hcrossing}\rput[l](5.1,1){$z$}
      \rput(6.5,.5){\flecheG}
      \rput(6.5,1.5){\flecheG}      
    \end{pspicture} \\
      & = & 
    \psset{unit=.5cm}
    \begin{pspicture}[.5](5.5,2)
      \rput[r](1,.5){$y$}\rput[r](1,1.5){$x$}
      \rput(1,.5){\flecheD}\rput(1,1.5){\flecheD}
      \psline(2,0)(2,2)\psline(3,0)(3,2)\psline(4.5,0)(4.5,2)
      \psline(2,.5)(3.5,.5)\psline(2,1.5)(3.5,1.5)
      \psline(4,.5)(4.5,.5)\psline(4,1.5)(4.5,1.5)
      \rput(3.75,1){$\dots$}
      \rput(5.5,.5){\flecheG}
      \rput(5.5,1.5){\flecheG}      
    \end{pspicture} \sigma(a\overline{z})
  \end{eqnarray*}
\endproof

The same method, together with the easy transformation
  \begin{equation}
    \psset{unit=.5cm}
    \begin{pspicture}[.45](2.5,1)
      \Hcrossing
      \psarc(1,.5){.5}{270}{90}
      \rput[r](.3,.5){$z$}
    \end{pspicture} = \left(\sigma(az)+\sigma(a^2)\right)
    \left(
    \begin{pspicture}[.45](1,1)
      \rput(0,1){\flecheD}\rput(1,0){\flecheG}
    \end{pspicture}\ +\ 
    \begin{pspicture}[.45](1,1)
      \rput(0,0){\flecheD}\rput(1,1){\flecheG}
    \end{pspicture}\right)
    \label{eq:boucle}
  \end{equation}
gives the following lemma.

\begin{lemm}\label{lemm:echange_boucle}
  \begin{eqnarray}
    \psset{unit=.5cm}
    \begin{pspicture}[.5](5.7,2)
      \rput[r](1,.5){$x$}\rput[r](1,1.5){$y$}
      \rput(1,.5){\flecheD}\rput(1,1.5){\flecheD}
      \psline(2,0)(2,2)\psline(3,0)(3,2)\psline(4.5,0)(4.5,2)
      \psline(2,.5)(3.5,.5)\psline(2,1.5)(3.5,1.5)
      \psline(4,.5)(4.5,.5)\psline(4,1.5)(4.5,1.5)
      \psarc(4.5,1){.5}{270}{90}
      \psline[linestyle=dotted](4.7,1)(5.7,1)
      \rput(3.75,1){$\dots$}
    \end{pspicture} &=&
    \psset{unit=.5cm}
    \frac{\sigma(a^2)+\sigma(x\overline{y})}{\sigma(a^2y\overline{x})}
    \begin{pspicture}[.5](5.7,2)
      \rput[r](1,.5){$y$}\rput[r](1,1.5){$x$}
      \rput(1,.5){\flecheD}\rput(1,1.5){\flecheD}
      \psline(2,0)(2,2)\psline(3,0)(3,2)\psline(4.5,0)(4.5,2)
      \psline(2,.5)(3.5,.5)\psline(2,1.5)(3.5,1.5)
      \psline(4,.5)(4.5,.5)\psline(4,1.5)(4.5,1.5)
      \psarc(4.5,1){.5}{270}{90}
      \psline[linestyle=dotted](4.7,1)(5.7,1)
      \rput(3.75,1){$\dots$}
    \end{pspicture} \label{eq:echange_boucle}\\
    \psset{unit=.5cm}
    \begin{pspicture}[.5](5.7,2)
      \rput[r](1,.5){$x$}\rput[r](1,1.5){$y$}
      \rput(1,.5){\flecheD}\rput(1,1.5){\flecheD}
      \psline(2,0)(2,2)\psline(3,0)(3,2)\psline(4.5,0)(4.5,2)
      \psline(2,.5)(3.5,.5)\psline(2,1.5)(3.5,1.5)
      \psline(4,.5)(4.5,.5)\psline(4,1.5)(4.5,1.5)
      \rput(4.5,1.5){\flecheD}
      \rput(5.5,.5){\flecheG}
      \rput(3.75,1){$\dots$}
    \end{pspicture} &=&
    \psset{unit=.5cm}
\frac{\sigma(x\overline{y})}{\sigma(a^2y\overline{x})}
    \begin{pspicture}[.5](5.7,2)
      \rput[r](1,.5){$y$}\rput[r](1,1.5){$x$}
      \rput(1,.5){\flecheD}\rput(1,1.5){\flecheD}
      \psline(2,0)(2,2)\psline(3,0)(3,2)\psline(4.5,0)(4.5,2)
      \psline(2,.5)(3.5,.5)\psline(2,1.5)(3.5,1.5)
      \psline(4,.5)(4.5,.5)\psline(4,1.5)(4.5,1.5)
      \rput(5.5,1.5){\flecheG}
      \rput(4.5,.5){\flecheD}
      \rput(3.75,1){$\dots$}
    \end{pspicture}  + 
\frac{\sigma(a^2)}{\sigma(a^2y\overline{x})}
    \begin{pspicture}[.5](5.7,2)
      \rput[r](1,.5){$y$}\rput[r](1,1.5){$x$}
      \rput(1,.5){\flecheD}\rput(1,1.5){\flecheD}
      \psline(2,0)(2,2)\psline(3,0)(3,2)\psline(4.5,0)(4.5,2)
      \psline(2,.5)(3.5,.5)\psline(2,1.5)(3.5,1.5)
      \psline(4,.5)(4.5,.5)\psline(4,1.5)(4.5,1.5)
      \rput(4.5,1.5){\flecheD}
      \rput(5.5,.5){\flecheG}
      \rput(3.75,1){$\dots$}
    \end{pspicture} \label{eq:echange_boucle_a}\\
    \psset{unit=.5cm}
    \begin{pspicture}[.5](5.7,2)
      \rput[r](1,.5){$x$}\rput[r](1,1.5){$y$}
      \rput(1,.5){\flecheD}\rput(1,1.5){\flecheD}
      \psline(2,0)(2,2)\psline(3,0)(3,2)\psline(4.5,0)(4.5,2)
      \psline(2,.5)(3.5,.5)\psline(2,1.5)(3.5,1.5)
      \psline(4,.5)(4.5,.5)\psline(4,1.5)(4.5,1.5)
      \rput(4.5,.5){\flecheD}
      \rput(5.5,1.5){\flecheG}
      \rput(3.75,1){$\dots$}
    \end{pspicture} &=&
    \psset{unit=.5cm}
\frac{\sigma(x\overline{y})}{\sigma(a^2y\overline{x})}
    \begin{pspicture}[.5](5.7,2)
      \rput[r](1,.5){$y$}\rput[r](1,1.5){$x$}
      \rput(1,.5){\flecheD}\rput(1,1.5){\flecheD}
      \psline(2,0)(2,2)\psline(3,0)(3,2)\psline(4.5,0)(4.5,2)
      \psline(2,.5)(3.5,.5)\psline(2,1.5)(3.5,1.5)
      \psline(4,.5)(4.5,.5)\psline(4,1.5)(4.5,1.5)
      \rput(5.5,.5){\flecheG}
      \rput(4.5,1.5){\flecheD}
      \rput(3.75,1){$\dots$}
    \end{pspicture}  + 
\frac{\sigma(a^2)}{\sigma(a^2y\overline{x})}
    \begin{pspicture}[.5](5.7,2)
      \rput[r](1,.5){$y$}\rput[r](1,1.5){$x$}
      \rput(1,.5){\flecheD}\rput(1,1.5){\flecheD}
      \psline(2,0)(2,2)\psline(3,0)(3,2)\psline(4.5,0)(4.5,2)
      \psline(2,.5)(3.5,.5)\psline(2,1.5)(3.5,1.5)
      \psline(4,.5)(4.5,.5)\psline(4,1.5)(4.5,1.5)
      \rput(4.5,.5){\flecheD}
      \rput(5.5,1.5){\flecheG}
      \rput(3.75,1){$\dots$}
    \end{pspicture} \label{eq:echange_boucle_b}
  \end{eqnarray}
\end{lemm}

We use Lemmas \ref{lemm:echange_lignes} and \ref{lemm:echange_boucle} to obtain symmetry properties of the partition functions, that we summarize below, where $m$ denotes either $2N$ or $2N+1$. 

\begin{lemm}\label{lemm:sym}
  The functions $Z(N;X_{2N})$ and 
  $Z_{\textsc{HT}}(2N+1;X_{2N},x,y)$ are symmetric separately in the two sets of variables $\{x_i,\ i\le N\}$ and $\{x_i,\ i\ge N+1\}$, the function 
  $Z_{\textsc{HT}}(2N;X_{2N-1},x,y)$ is symmetric separately in the two sets of variables $\{x_i,\ i\le N-1\}$ and $\{x_i,\ i\ge N\}$, 
  and the functions 
  $Z_{\textsc{QT}}(2m;X_{N-1},x,y)$ are symmetric in their variables $x_i$.

Moreover, $Z_{\textsc{QT}}(4N+2;\dots)$ is symmetric in its variables $x$ and
$y$, and we have a pseudo-symmetry for $Z_{\textsc{QT}}(4N;\dots)$ and $Z_{\textsc{HT}}(2N;\dots)$:
  \begin{eqnarray}
    Z_{\textsc{QT}}(4N;X_{2N-1},x,y) & = & 
    \frac{\sigma(a^2)+\sigma(x\overline{y})}{\sigma(a^2y\overline{x})}
    Z_{\textsc{QT}}(4N;X_{2N-1},y,x),
    \label{eq:sym_QT4N_xy}\\
    Z_{\textsc{HT}}(2N;X_{2N-1},x,y) & = & 
    \frac{\sigma(a^2)+\sigma(x\overline{y})}{\sigma(a^2y\overline{x})}
    Z_{\textsc{HT}}(2N;X_{2N-1},y,x).
    \label{eq:sym_HT2N_xy}
  \end{eqnarray}
\end{lemm}

\proof
  For $Z(N;\dots)$ and $Z_{\textsc{HT}}(m;\dots)$, the symmetry in two ``consecutive'' 
  variables $x_i$ and $x_{i+1}$ is a direct consequence of
  Lemma~\ref{lemm:echange_lignes}. For $Z_{\textsc{QT}}(2m;\dots)$, we again apply 
  Lemma~\ref{lemm:echange_lignes} together with the easy observations:
\begin{equation}\label{eq:2dots}
  \psset{unit=5mm}
    \begin{pspicture}[.5](2,2)
      \psline(1,0)(1,2)\psline(0,1)(2,1)
    \end{pspicture}
    \ =\  
     \begin{pspicture}[.5](2,2)
       \psline(1,0)(1,2)\psline(0,1)(2,1)
       \psdot[dotsize=.2](.5,1)
       \psdot[dotsize=.2](1,1.5)
       \psdot[dotsize=.2](1,.5)
       \psdot[dotsize=.2](1.5,1)
     \end{pspicture} 
\ \ \ \ {\rm and}\ \ \ \ 
     \begin{pspicture}[.5](2,2)
       \psline(0,1)(2,1)
     \end{pspicture}
    \ =\  
     \begin{pspicture}[.5](2,2)
       \psline(0,1)(2,1)
       \psdot[dotsize=.2](.5,1)
       \psdot[dotsize=.2](1.5,1)
     \end{pspicture}
\end{equation}
which allow us to bring the Yang-Baxter triangle through the divalent vertices of Figure \ref{fig:ZQT}.

For the (pseudo-)symmetries in $(x,y)$, let us deal with $Z_{\textsc{QT}}(4N;\dots)$, the other cases being 
similar or simpler. We use equation \pref{eq:2dots} to put together the lines of parameter $x$ and $y$:

\begin{eqnarray*}
Z_{\textsc{QT}}(4N;X_{2N-1},x,y) & = & 
\psset{unit=.3cm}\begin{pspicture}[.4](10,13)
\rput(0,0){\QTIce{6}}
\rput[r](0,6){$x$}\rput[t](6,-.1){$y$}
\end{pspicture}\\
& = &
\psset{unit=.3cm}
  \begin{pspicture}[.4](12,12)
    \rput(1,1){\IceGrid{7}{5}}
    \rput(0,1){\colD{5}}
    \rput(1,0){\linB{7}}
    \rput[t](6,-.1){$y$}
    \rput[t](7,-.1){$x$}
    \psarc(6.5,5){.5}{0}{180}
    \psline[linestyle=dotted](6.5,5.2)(6.5,6.2)
    \multido{\i=1+1}{5}{
      \psline(\i,5)(\i,6)
      \psarc(6,6){\i}{90}{180}
      \psarc(7,6){\i}{270}{90}
      }
      \multido{\i=7+1}{5}{\psline(6,\i)(7,\i)}
    \rput(7,6){\SpecialCoor\multido{\i=1+1}{5}{\psdots[dotstyle=*](\i;45)}}
    \end{pspicture}
\end{eqnarray*}
and then apply Lemma \ref{lemm:echange_boucle}.
\endproof

It should be clear that we  have analogous properties for the even and odd parts of the partition functions.
The next (and last) symmetry property, proved by Stroganov \cite{IKdet}, appears when the parameter $a$ equals the special value $\omega_6=\exp(i\pi/3)$. An elementary proof of this result has recently been given in \cite{aval}.

\begin{lemm}\label{lemm:symZomega6}
When $a=\omega_6$, the partition function $Z(N;X_{2N})$ is symmetric in {\bf all} its variables.
\end{lemm}

\subsection{Specializations, recurrences}

The aim of this section is to give the value of the partition functions in some specializations of the variable $x$ or $y$. The first result is due to Kuperberg, the others are very similar.

\begin{lemm}{\em [specialization of $Z$; Kuperberg]}
\label{lemm:specZ}
If we denote
\begin{eqnarray*}
  A(x_{N+1},X_{2N}\backslash \{x_1,x_{N+1}\}) & = & 
    \prod_{2\leq k\leq N} \sigma(a x_k \overline{x}_{N+1})
    \prod_{N+1\leq k\leq 2N} \sigma(a^2 x_{N+1} \overline{x}_k),\\
  \overline{A}(x_{N+1},X_{2N}\backslash \{x_1,x_{N+1}\}) & = & 
    \prod_{2\leq k\leq N} \sigma(a x_{N+1} \overline{x}_k) 
    \prod_{N+1\leq k\leq 2N} \sigma(a^2 x_k \overline{x}_{N+1}),
\end{eqnarray*}
then we have:
  \begin{eqnarray}
 Z(N;{\bf \overline{a}x_{N+1}},X_{2N}\backslash x_1) & = & 
\overline{A}(x_{N+1},X_{2N}\backslash \{x_1,x_{N+1}\}) Z(N-1;X_{2N}\backslash \{x_1,x_{N+1}\}), \label{eq:Zbax}\\
 Z(N;{\bf ax_{N+1}},X_{2N}\backslash x_1) & = & 
A(x_{N+1},X_{2N}\backslash \{x_1,x_{N+1}\}) 
 Z(N-1; X_{2N}\backslash \{x_1,x_{N+1}\}). \label{eq:Zax}
  \end{eqnarray}
\end{lemm}

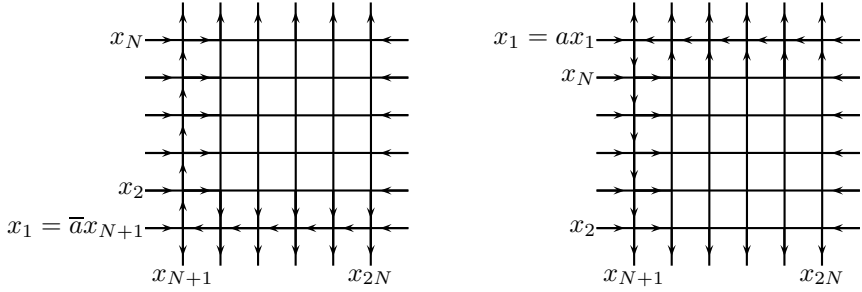
\begin{figure}[htbp]
    \begin{center}\psset{unit=.5cm}
      \begin{pspicture}(-1,-1)(19,8)
	\rput(0,0){\IceSquare{6}}
	\rput[r](0,1){$x_1=\overline{a}x_{N+1}$}\rput[r](0,6){$x_{N}$}
        \rput[r](0,2){$x_2$}
        \rput[t](1,-.1){$x_{N+1}$}\rput[t](6,-.1){$x_{2N}$}
	\multido{\i=1+1}{5}{\rput(1,\i){\flecheH}}
        \multido{\i=2+1}{5}{\rput(1,\i){\flecheD}\rput(\i,1){\flecheG}\rput(\i,2){\flecheB}}

	\rput(12,0){%
	  \IceSquare{6}
          \rput[r](0,6){$x_1=ax_{1}$}\rput[r](0,5){$x_{N}$}\rput[r](0,1){$x_2$}
          \rput[t](1,-.1){$x_{N+1}$}\rput[t](6,-.1){$x_{2N}$}
          \multido{\i=2+1}{5}{\rput(1,\i){\flecheB}\rput(\i,5){\flecheH}
             \rput(\i,6){\flecheG}}
          \multido{\i=1+1}{5}{\rput(1,\i){\flecheD}}
	}
      \end{pspicture}
      \caption{Fixed edges for (\ref{eq:Zbax}) on the left and
      (\ref{eq:Zax}) on the right}
    \label{fig:recZ}
    \end{center}
  \end{figure}

\proof
We recall the method to prove equation \pref{eq:Zbax}. We observe that when $x_1=\bar a x_{N+1}$, the parameter of the vertex at the crossing of the two lines of parameter $x_1$ and $x_{N+1}$ is $\bar a$. Thus the weight of this vertex is $\si(a\ba)=\si(1)=0$ unless the orientation of this vertex is the second on Figure \ref{fig:poids_6V}. But this orientation implies the orientation of all vertices in the row $x_1$ and in the column $x_{N+1}$, as shown on Figure \ref{fig:recZ}. 
The non-fixed part gives the partition function $Z$ in size $N-1$, without parameters $x_1$ and
  $x_{N+1}$, and the weights of the fixed part gives the factor $\overline{A}(\dots)$.

The case of \pref{eq:Zax} is similar, after using Lemma \ref{lemm:sym} to put the line $x_{N+1}$ at the top of the grid.

\endproof

We will need the following application of the Yang-Baxter equation, which allows, under certain condition, a line with a change of parameter to go through a grid.

\begin{lemm}\label{lemm:tour_grille}
  \begin{equation}
    \psset{unit=.3cm}
    \begin{pspicture}[.5](9,8)
      \multido{\i=2+1}{6}{\psline(\i,0)(\i,8)}
      \multido{\i=2+1}{5}{\psline(0,\i)(9,\i)}
      \psline(1,7)(7,7)\psline(8,1)(8,6)
      \psarc(7,6){1}{0}{90}
      \psline[linestyle=dotted](7.2,6.2)(8.5,7.5)
      \rput[r](1,7){$x$}\rput[t](8,.9){$\overline{a}x$}
    \end{pspicture}
    = 
    \begin{pspicture}[.5](9,8)
      \multido{\i=2+1}{6}{\psline(\i,0)(\i,8)}
      \multido{\i=2+1}{5}{\psline(0,\i)(9,\i)}
      \psline(1,7)(1,2)\psline(2,1)(8,1)
      \psarc(2,2){1}{180}{270}
      \psline[linestyle=dotted](1.8,1.8)(.5,.5)
      \rput[l](8,1){$x$}\rput[b](1,7){$\overline{a}x$}
    \end{pspicture}
  \label{eq:tour_grille}
  \end{equation}
\end{lemm}

\proof
 We iteratively apply Lemma \ref{lemm:YB} on the rows, and row by row:
  \begin{eqnarray*}
    \psset{unit=.3cm}
    \begin{pspicture}[.5](10,4)
      \psline(0,2)(10,2)
      \multido{\i=2+2}{4}{\psline(\i,0)(\i,4)}
      \psline(1,3)(8,3)\psline(9,2)(9,1)\psarc(8,2){1}{0}{90}
      \psline[linestyle=dotted](8.2,2.2)(9.5,3.5)
      \rput[r](1,3){$x$}\rput[t](9,.9){$\overline{a}x$}
    \end{pspicture}
    & = & 
    \psset{unit=.3cm}
    \begin{pspicture}[.5](10,4)
      \psline(0,2)(10,2)
      \multido{\i=2+2}{4}{\psline(\i,0)(\i,4)}
      \psline(1,3)(6,3)\psline(8,1)(9,1)\psarc(6,2){1}{0}{90}
      \psarc(8,2){1}{180}{270}
      \psline[linestyle=dotted](6.2,2.2)(7.5,3.5)
      \psline[linestyle=dotted](7.8,1.8)(6.5,.5)
      \rput[r](1,3){$x$}
      \rput[l](9,1){$x$}
      \rput[r](6.8,1.5){$\overline{a}x$}
    \end{pspicture}\\
    & = & 
    \psset{unit=.3cm}
    \begin{pspicture}[.5](10,4)
      \psline(0,2)(10,2)
      \multido{\i=2+2}{4}{\psline(\i,0)(\i,4)}
      \psline(1,3)(4,3)\psline(6,1)(9,1)\psarc(4,2){1}{0}{90}
      \psarc(6,2){1}{180}{270}
      \psline[linestyle=dotted](4.2,2.2)(5.5,3.5)
      \psline[linestyle=dotted](5.8,1.8)(4.5,.5)
      \rput[r](1,3){$x$}
      \rput[l](9,1){$x$}
      \rput[r](4.8,1.5){$\overline{a}x$}
    \end{pspicture}\\
    & = & 
    \psset{unit=.3cm}
    \begin{pspicture}[.5](10,4)
      \psline(0,2)(10,2)
      \multido{\i=2+2}{4}{\psline(\i,0)(\i,4)}
      \psline(1,3)(4,3)\psline(6,1)(9,1)\psarc(4,2){1}{0}{90}
      \psarc(6,2){1}{180}{270}
      \psline[linestyle=dotted](4.2,2.2)(5.5,3.5)
      \psline[linestyle=dotted](5.8,1.8)(4.5,.5)
      \rput[r](1,3){$x$}
      \rput[l](9,1){$x$}
      \rput[r](4.8,1.5){$\overline{a}x$}
    \end{pspicture}\\
    & = & 
    \psset{unit=.3cm}
    \begin{pspicture}[.5](10,4)
      \psline(0,2)(10,2)
      \multido{\i=2+2}{4}{\psline(\i,0)(\i,4)}
      \psline(1,3)(1,2)\psline(2,1)(9,1)\psarc(2,2){1}{180}{270}
      \psline[linestyle=dotted](1.8,1.8)(.5,.5)
      \rput[b](1,3){$\overline{a}x$}\rput[l](9,1){$x$}
    \end{pspicture}.
  \end{eqnarray*}

\endproof

\begin{lemm}{\em [specialization of $Z_{\textsc{HT}}$]}
\label{lemm:specZHT}
If we denote
 \begin{eqnarray*}
   A_{H}^{1}(x_1,X_{2N}\backslash x_1) & = &    
     \prod_{1\leq k\leq N} \sigma(a^2 x_1 \overline{x}_k)
     \prod_{N+1\leq k\leq 2N} \sigma(a x_k\overline{x}_1),\\
   \overline{A}_{H}^{1}(x_1,X_{2N}\backslash x_1) & = & 
     \prod_{1\leq k\leq N} \sigma(a^2x_k\overline{x}_1)
     \prod_{N+1\leq k\leq 2N} \sigma(ax_1\overline{x}_k),\\
   A_{H}^{0}(x_{N},X_{2N-1}\backslash x_{N}) & = & 
     \prod_{1\leq k\leq N-1}\sigma(a x_k \overline{x}_N)
     \prod_{N\leq k\leq 2N-1} \sigma(a^2 x_N \overline{x}_k),\\
   \overline{A}_{H}^{0}(x_{N},X_{2N-1}\backslash x_{N}) & = & 
     \prod_{1\leq k\leq N-1} \sigma(ax_N\overline{x}_k)
     \prod_{N\leq k\leq 2N-1} \sigma(a^2x_k\overline{x}_N),
 \end{eqnarray*}
  then for $\star=\downright,\upleft,\upc,\downc$ and
  $\square=\downc,\upc,\downright,\upleft$ respectively, we have
  \begin{eqnarray}
\!\!\!\!\!\!\!\!\!\!\!\!Z_{\textsc{HT}}^{\star}(2N+1;X_{2N},x,\mathbf{ax_1}) &\!\!\!\!=\!\!\!\!&
    A_{H}^{1}(x_1,X_{2N}\backslash x_{1}) 
    Z_{\textsc{HT}}^{\square}(2N;X_{2N}\backslash x_1, x_1,x),
    \label{equa:Zht_impair_ax}\\
\!\!\!\!\!\!\!\!\!\!\!\!Z_{\textsc{HT}}^{\square}(2N+1;X_{2N},x,{\bf \overline{a}x_1}) &\!\!\!\!=\!\!\!\!&
    \overline{A}_{H}^{1}(x_1,X_{2N}\backslash x_1)
     Z_{\textsc{HT}}^{\star}(2N;X_{2N}\backslash x_1,x,x_1),
    \label{equa:Zht_impair_bax}\\
\!\!\!\!\!\!\!\!\!\!\!\!Z_{\textsc{HT}}^{\star}(2N;X_{2N-1},x,{\bf ax_N}) &\!\!\!\!=\!\!\!\!&
  \sigma(ax\overline{x}_{N}) A_{H}^{0}(x_{N},X_{2N-1}\backslash x_{N})
  Z_{\textsc{HT}}^{\square}(2N-1;X_{2N-1}\backslash x_N,x,x_N),   
    \label{equa:Zht_pair_ax}\\
\!\!\!\!\!\!\!\!\!\!\!\!Z_{\textsc{HT}}^{\square}(2N;X_{2N-1},{\bf \overline{a}x_N},y) &\!\!\!\!=\!\!\!\!&
  \sigma(ax_N\overline{y}) 
  \overline{A}_{H}^{0}(x_{N},X_{2N-1}\backslash x_{N})
   Z_{\textsc{HT}}^{\star}(2N-1;X_{2N-1}\backslash x_N,y,x_N).
    \label{equa:Zht_pair_bax}
  \end{eqnarray}
\end{lemm}

\proof

The proof is similar to the previous one, with the difference 
that before looking at fixed edges, we need to multiply the partition function by a given factor; 
we interpret this operation by a modification of the graph of the ice model, and apply 
Lemma~\ref{lemm:tour_grille}. It turns out that in each case, the additional factors 
are exactly cancelled by the weights of fixed vertices.

  To prove (\ref{equa:Zht_impair_ax}), we multiply the left-hand side by
  \begin{displaymath}
    \prod_{N+1\leq k\leq 2N} \sigma(a^2 x_k \overline{y}),
  \end{displaymath}
  which is equivalent to adding to the line of parameter $y$ a new line
  $\overline{a}y$ just below the grid;
  Lemma~\ref{lemm:tour_grille} transforms the graph of
  Figure~\ref{fig:Zht_impair_ax}(a) into the graph of Figure \ref{fig:Zht_impair_ax}(b). 
  When we put  $y=ax_1$, we get the indicated fixed edges, which gives as partition function
  \begin{displaymath}
    \prod_{N+1\leq k\leq 2N} \sigma^2(ax_k\overline{x}_1)
    \prod_{1\leq k\leq N} \sigma(a^2x_1\overline{x}_k)
    Z_{\textsc{HT}}(2N;X_{2N}\backslash x_1, x_1,x).
  \end{displaymath}

  \begin{figure}[htbp]
    \begin{center}
      \psset{unit=5mm}
      \begin{pspicture}[.5](-1,-1)(9,11)
	\rput(1,1){\IceGrid{4}{10}}
	\rput(1,0){\linB{4}}
	\rput(1,1){\linB{4}}
	\rput(0,2){\colD{9}}
	\rput(1,10){\linH{4}}
	\psline(5,2)(5,5)
	\multido{\i=1+1}{4}{\rput(\i,1){\flecheG}}
	\multido{\i=2+1}{4}{\psline(4,\i)(5,\i)}
	\multido{\i=7+1}{4}{\psline(4,\i)(5,\i)}
	\multido{\i=1+1}{4}{\psarc(5,6){\i}{270}{90}}
	\psarc(4,2){1}{270}{0}
	\psarc(4,5){1}{0}{90}
	\SpecialCoor
	\rput(4,5){\psline[linestyle=dotted](.5;45)(1.5;45)}
	\rput(4,2){\psline[linestyle=dotted](.5;315)(1.5;315)}
	\rput[r](0,2){$x_1$}
	\rput[r](0,5){$x_{N}$}
	\rput[r](0,6){$x$}
	\rput[r](0,1){$\overline{a}y$}
	\rput[t](1,-.1){$x_{N+1}$}
	\rput[t](4,-.1){$x_{2N}$}
	\rput[l](5.1,1.2){$y$}
	\rput[c](4.5,-1.5){(a)}
      \end{pspicture} \ \ 
      \begin{pspicture}[.5](-1,-1)(10,11)
	\rput(2,1){\IceGrid{4}{10}}
	\rput(1,0){\linB{5}}
	\rput(2,1){\linB{4}}
	\rput(2,9){\linH{4}}
	\rput(2,10){\linH{4}}
	\rput(0,1){\colD{4}}
	\rput(1,2){\colD{3}}
	\rput(1,6){\colD{5}}
	\multido{\n=.5+1.0}{5}{\psarc(5,5.5){\n}{270}{90}}
	\multido{\i=1+1}{3}{\rput(1,\i){\flecheH}}
	\multido{\i=2+1}{4}{\rput(\i,1){\flecheG}}
	\multido{\i=2+1}{3}{\rput(\i,10){\flecheD}}
	\psarc(2,4){1}{90}{180}
	\psarc{<-}(5,5.5){4.5}{0}{90}
	\SpecialCoor
	\rput(2,4){\psline[linestyle=dotted](.5;135)(1.5;135)}
	\rput(5,5.5){\psline[linestyle=dotted](.2;0)(1;2)}
	\rput[r](0,1){$x_1$}
	\rput[r](0,4){$x_N$}
	\rput[r](1,6){$x$}
	\rput[t](1,-.1){$y$}
	\rput[t](2.5,-.1){$x_{N+1}$}
	\rput[t](5,-.1){$x_{2N}$}
	\rput[b](.8,5){$\overline{a}y=x_1$}
	\rput[c](5,-1.5){(b)}
      \end{pspicture}
    \end{center}
    \caption{Proof of (\ref{equa:Zht_impair_ax})}\label{fig:Zht_impair_ax}
  \end{figure}

  Since
  $a^2x_k\overline{y}=ax_k\overline{x}_1$, the equation simplifies. 
  To conclude, we observe that if we start with an edge going out from the crossing $x/x_{2N}$ 
  (function $Z_{\textsc{HT}}^{\downright}$) we get at the end the same orientation (function
  $Z_{\textsc{HT}}^{\downc}$).
\endproof

\begin{lemm}{\em [specialization of $Z_{\textsc{QT}}$]}
\label{lemm:specZQT}
If we denote
\begin{eqnarray*}
  \overline{A}_{Q}(x_1,X_{m-1}\backslash x_1) & = & 
  \prod_{1\leq k\leq m-1} \sigma(a^2 x_k\overline{x}_1) 
   \sigma(ax_1 \overline{x}_k),\\
  A_{Q}(x_1;X_{m-1}\backslash x_1) & = & 
  \prod_{1\leq k\leq m-1} \sigma(a^2 x_1 \overline{x}_k)
\sigma(a x_k \overline{x}_1),
\end{eqnarray*}
then for $\star = \convcorner,\divcorner,\downright,\upleft$ and
$\square=\upleft,\downright,\convcorner,\divcorner$ respectively, we have:
\begin{eqnarray}
\!\!\!\!\!\!  Z_{\textsc{QT}}^{\star}(2m;X_{m-1},\mathbf{\overline{a}x_1},y) &=&
  \sigma(ax_1\overline{y}) \overline{A}_{Q}(x_1,X_{m-1})
  Z_{\textsc{QT}}^{\square}(2m-2;X_{m-1}\backslash x_1, y,x_1),
  \label{eq:Zqt_spe_bax}\\
\!\!\!\!\!\!   Z_{\textsc{QT}}^{\square}(2m;X_{m-1},x,\mathbf{ax_1}) & = & 
\sigma(ax\overline{x}_1) A_{Q}(x_1;X_{m-1}\backslash x_1)
  Z_{\textsc{QT}}^{\star}(2m-2;X_{m-1}\backslash x_1, x_1,x).
  \label{eq:Zqt_spe_ax}
\end{eqnarray}
\end{lemm}

\proof
Similar to the proof of Lemma \ref{lemm:specZHT}.
\endproof

\begin{rema}\label{rema:sym}
By using the (pseudo-)symmetry in $(x,y)$, we may transform any specialization of the variable $y$ into a specialization of the variable $x$. Moreover, by using Lemma \ref{lemm:sym} and (when $a=\omega_6$) Lemma \ref{lemm:symZomega6}, we obtain for $Z$, $Z_{\textsc{HT}}$  and $Z_{\textsc{QT}}$, $2N$ specializations. Now we have to compare them.

\end{rema}

\subsection{Special value of the parameter $a$; conclusion}

When $a=\omega_6=\exp(i\pi/3)$, two new ingredients may be used. The first one is Lemma \ref{lemm:symZomega6}, as mentioned in Remark \ref{rema:sym}. The second one is that with this special value of $a$:
\begin{equation}\label{spec}
\si(a)=\si(a^2) \ \ \ \ \ \ \ \si(a^2x)=-\si(\ba x)=\si(a\bx).
\end{equation}
which implies that the products appearing in Lemmas \ref{lemm:specZ}, \ref{lemm:specZHT} and \ref{lemm:specZQT} may be written in a more compact way:
\begin{eqnarray*}
  A(x_{N+1},X_{2N}\backslash\{x_1,x_{N+1}\}) & = & \sigma(a)
  \prod_{k\neq 1,N+1} \sigma(a x_k \overline{x}_{N+1}),\\
  \overline{A}(x_{N+1},X_{2N}\backslash \{x_1,x_{N+1}\}) & = & 
  \sigma(a) \prod_{k\neq 1,N+1} \sigma(a x_{N+1} \overline{x}_k),\\
  A_{H}^{1}(x_1,X_{2N}\backslash x_1) & = & 
  \prod_{1\leq k\leq 2N} \sigma(a x_{k} \overline{x}_1),\\
  \overline{A}_{H}^{1}(x_1,X_{2N}\backslash x_1) & = & 
  \prod_{1\leq k\leq 2N} \sigma(a x_1 \overline{x}_k),\\
  A_{H}^{0}(x_N,X_{2N-1}\backslash x_N) & = & 
  \prod_{1\leq k\leq 2N-1} \sigma(a x_{k} \overline{x}_{N}),\\
  \overline{A}_{H}^{0}(x_N,X_{2N-1}\backslash x_{N}) &=&
  \prod_{1\leq k\leq 2N-1} \sigma(a x_{N} \overline{x}_k),\\
  \overline{A}_{Q}(x_1,X_{m-1}\backslash x_1) & = &
  \prod_{1\leq k\leq m-1} \sigma^2(a x_1 \overline{x}_k),\\
  A_{Q}(x_1,X_{m-1}\backslash x_1) & = & 
  \prod_{1\leq k\leq m-1} \sigma^2(a x_k \overline{x}_1).
\end{eqnarray*}

Thus we get by comparing:
\begin{eqnarray*}
  A(x_i, X_{2N}\backslash x_i, x)
  A_{H}^{1}(x_i,X_{2N}\backslash x_i)& = &\sigma(a x
  \overline{x}_i) A_{Q}(x_i,X_{2N}\backslash x_i)\\
  \overline{A}(x_i,X_{2N}\backslash x_i,x)
  \overline{A}_{H}^{1}(x_i,X_{2N}\backslash x_i) & = &
  \sigma(a x_i \overline{x})
  \overline{A}_{Q}(x_i,X_{2N}\backslash x_i),
\end{eqnarray*}
whence (\ref{eq:main1.1}) and (\ref{eq:main1.2}) imply that
(\ref{eq:main2.1}) and (\ref{eq:main2.2}) are true (in size
$4N+2$) for the $2N$ specializations $x=a^{\pm 1}x_i$ ($1\leq i\leq
N$). It is enough to prove (\ref{eq:main2.2}) (Laurent polynomials of half-width
 $2N-1$), but we still need one specialization to get
 (\ref{eq:main2.1}) (half-width $2N$).

For (\ref{eq:main1.1}) and (\ref{eq:main1.2}), we observe the same kind of simplification 
\begin{displaymath}
  A(x_i, X_{2N-1}\backslash x_i) \sigma(a x \overline{x}_i)
  A_{H}^{0}(x_i,X_{2N-1}\backslash x_i) = \sigma(a x \overline{x}_i)
  A_{Q}(x_i,X_{2N-1}\backslash x_i),
\end{displaymath}
whence (\ref{eq:main2.2}) and (\ref{eq:main2.1}) for the size $4N-2$
imply that (\ref{eq:main1.1}) and (\ref{eq:main1.2}) are true for the 
 $N$ specializations $x=ax_i$, $N\leq i\leq 2N-1$. We obtain in the same way
the coincidence for the $N$ specializations $x=\overline{a}x_i$, $N\leq i\leq
2N-1$. Thus we have $2N$ specialiations of $x$: it is enough both for
 (\ref{eq:main1.1}) (half-width  $2N-1$), and for (\ref{eq:main1.2}) (half-width $2N-2$).

At this point, we have \emph{almost} proved
\begin{center}
  ((\ref{eq:main1.1}) and (\ref{eq:main1.2}), in size $4N$) $\Longrightarrow$
  ((\ref{eq:main2.1}) and (\ref{eq:main2.2}), in size $4N+2$) $\Longrightarrow$
  ((\ref{eq:main1.1}) and (\ref{eq:main1.2}), in size $4N+4$);
\end{center}
\emph{almost}, because we still need {\em one} specialization for
(\ref{eq:main2.1}).

We get this missing specialization, not directly for 
$\ZqDownright$, $\ZqUpleft$, $\ZhDownright$ and $\ZhUpleft$, but for the original series $Z_{\textsc{QT}}(4N+2;X_{2N},x,y)$ and
$Z_{\textsc{HT}}(2N+1;X_{2N},x,y)$: indeed if we set
$x=ay$ we may apply
Lemma~\ref{lemm:tour_grille}.

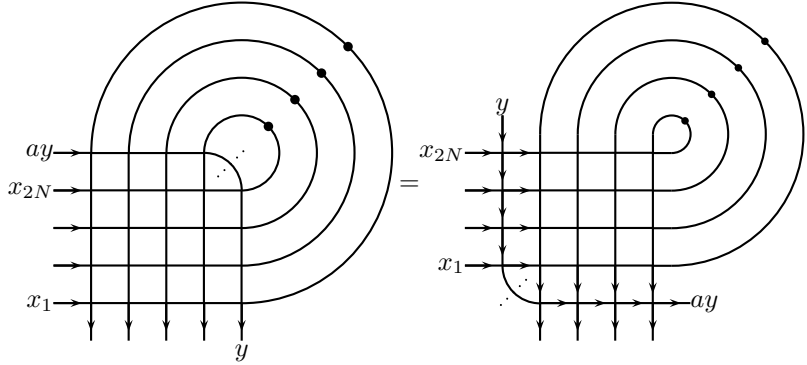
$$  \psset{unit=5mm}
  \begin{pspicture}[.5](-1,-1)(9,9)
    \rput(0,0){\qQTIce{5}}
    \rput[r](0,1){$x_1$}
    \rput[r](0,4){$x_{2N}$}
    \rput[r](0,5){$ay$}
    \rput[t](5,-.1){$y$}
  \end{pspicture}
 =  
  \psset{unit=5mm}
  \begin{pspicture}[.5](-1,-1)(9,9)
    \rput(2,2){\IceGrid{4}{4}}
    \rput(2,1){\linB{4}}
    \rput(2,0){\linB{4}}
    \rput(1,2){\colD{4}}
    \rput(0,2){\colD{4}}
    \multido{\i=2+1}{4}{\rput(\i,1){\flecheD}}
    \multido{\i=3+1}{4}{\rput(1,\i){\flecheB}}
    \multido{\i=2+1}{4}{\psline(\i,5)(\i,5.5)\psline(5,\i)(5.5,\i)}
    \multido{\n=.5+1.0}{4}{\psarc(5.5,5.5){\n}{270}{180}}
    \SpecialCoor
    \rput(2,2){\psline[linestyle=dotted](.5;225)(1.5;225)}
    \rput(5.5,5.5){\psdots[dotsize=1mm](.5;45)(1.5;45)(2.5;45)(3.5;45)}
    \psarc(2,2){1}{180}{270}
    \rput[b](1,6){$y$}
    \rput[r](0,2){$x_1$}
    \rput[r](0,5){$x_{2N}$}
    \rput[l](6,1){$ay$}
  \end{pspicture}$$
$$Z_{\textsc{QT}}(4N+2;X_{2N},{\bf ay},y)  =  
\sigma(a) \prod_{1\leq k\leq 2N} \sigma(a x_k \overline{y}) \sigma(a^2 y \overline{x}_k) Z_{\textsc{QT}}(4N;X_{2N}\backslash x_{2N}, x_{2N},x_{2N})$$
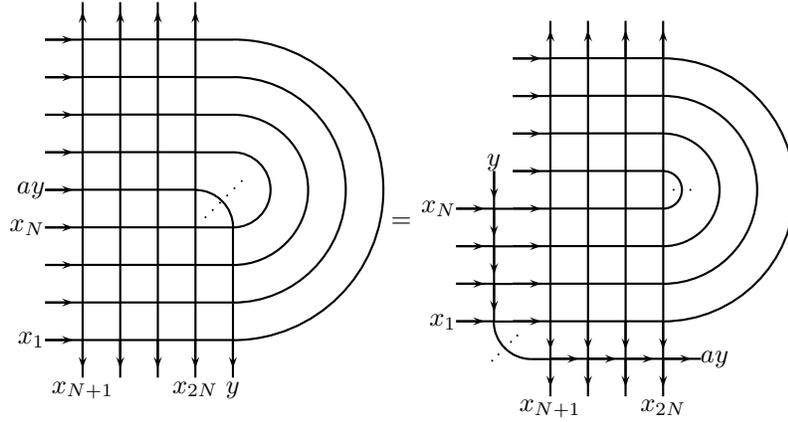
$$  \psset{unit=5mm}
  \begin{pspicture}[.5](-1,-1)(9,9)
    \rput(0,0){\HTIceOdd{4}}
    \rput[r](0,1){$x_1$}
    \rput[r](0,4){$x_N$}
    \rput[r](0,5){$ ay$}
    \rput[t](1,-.1){$x_{N+1}$}
    \rput[t](4,-.1){$x_{2N}$}
    \rput[t](5,-.1){${y}$}
  \end{pspicture}
  =  
  \psset{unit=5mm}
  \begin{pspicture}[.5](-.5,-1)(9,10)
    \rput(2,1){\HTIceEven{4}}
    \rput(3,0){\linB{4}}
    \rput(.5,2){\colD{4}}
    \multido{\i=2+1}{4}{\psline(1.5,\i)(2,\i)}
    \psline(2.5,1)(3,1)
    \multido{\i=3+1}{4}{\rput(\i,1){\flecheD}}
    \multido{\i=3+1}{4}{\rput(1.5,\i){\flecheB}}
    \SpecialCoor
    \psarc(2.5,2){1}{180}{270}
    \rput(2.5,2){\psline[linestyle=dotted](.5;225)(1.5;225)}
    \rput[r](.5,2){$x_1$}
    \rput[r](.5,5){$x_{N}$}
    \rput[b](1.5,6){$y$}
    \rput[t](3,-.1){$x_{N+1}$}
    \rput[t](6,-.1){$x_{2N}$}
    \rput[l](7,1){$ay$}
  \end{pspicture}$$
$$  Z_{\textsc{HT}}(2N+1;X_{2N},{\bf ay},y)  = 
\left(\prod_{1\leq k\leq N} \sigma(a x_k \overline{y}) 
\prod_{N+1\leq k\leq 2N} \sigma(a^2 y \overline{x}_k)\right)
  Z_{\textsc{HT}}(2N; X_{2N}\backslash x_{N}, x_N,x_N)
$$

This way, we get another point where
(\ref{eq:main2}) is true, and thus, because we already have
(\ref{eq:main2.2}), by difference we obtain that (\ref{eq:main2.1}) 
holds for $y=\overline{a}x$.

This completes the proof of Theorem \ref{theo:main}.


\end{document}